\documentclass[graybox]{svmult}
\usepackage[T1]{fontenc}
\usepackage{newtxtext,newtxmath}
\usepackage[bottom]{footmisc}

\usepackage{booktabs}
\usepackage{subcaption}
\usepackage{mathtools}
\mathtoolsset{centercolon,mathic}
\usepackage{tikz-cd}
\usetikzlibrary{positioning,decorations.markings}

\pgfdeclarearrow{
  name = txto,
  setup code = {
    \pgfarrowssettipend{1.9285\pgflinewidth}
    \pgfarrowssetbackend{-2.125\pgflinewidth}
    \pgfarrowssetlineend{0pt}
    \pgfarrowssetvisualbackend{0pt}
    \pgfarrowshullpoint{1.9285\pgflinewidth}{0pt}
    \pgfarrowsupperhullpoint{-1.5714\pgflinewidth}{3.6607\pgflinewidth}
    \pgfarrowsupperhullpoint{-2.125\pgflinewidth}{3.0714\pgflinewidth}
  },
  drawing code = {
    \pgfsetdash{}{0pt}%
    \pgfpathmoveto{\pgfpoint{1.9285\pgflinewidth}{0pt}}%
    \pgfpathlineto{\pgfpoint{-1.5714\pgflinewidth}{3.6607\pgflinewidth}}%
    \pgfpathlineto{\pgfpoint{-2.125\pgflinewidth}{3.0714\pgflinewidth}}%
    \pgfpathlineto{\pgfpoint{0pt}{0.5\pgflinewidth}}%
    \pgfpathlineto{\pgfpoint{0pt}{-0.5\pgflinewidth}}%
    \pgfpathlineto{\pgfpoint{-2.125\pgflinewidth}{-3.0714\pgflinewidth}}%
    \pgfpathlineto{\pgfpoint{-1.5714\pgflinewidth}{-3.6607\pgflinewidth}}%
    \pgfpathclose%
    \pgfusepathqfill
  }
}

\tikzset{>=txto}
\newlength\myfontsize
\newlength\mythinwidth
\newlength\mythickwidth
\newlength\mydoubledist
\tikzset{every picture/.style={line width={\mythinwidth}}}
\tikzset{thin/.style={line width={\mythinwidth}}}
\tikzset{thick/.style={line width={\mythickwidth}}}
\tikzset{arrow/.style={->,node font=\small}}
\tikzcdset{arrow style=tikz}

\usepackage[hidelinks]{hyperref}

\newcommand*{\Z}{\mathbb{Z}}

\newcommand*{\calC}{\mathcal{C}}
\newcommand*{\calD}{\mathcal{D}}
\newcommand*{\calS}{\mathcal{S}}
\newcommand*{\Ntwo}{\mathbb{N}_{-2}}

\newcommand*{\inftyone}{\ensuremath{(\infty,1)}}

\newcommand*{\SmoothInftyGpd}{\ensuremath{\mathbf{Smooth}\infty\mathbf{Gpd}}}
\newcommand*{\Type}{\ensuremath{\mathbf{Type}}}
\newcommand*{\sType}{\ensuremath{\mathbf{sType}}}
\newcommand*{\ssType}{\ensuremath{\mathbf{ssType}}}
\newcommand*{\Prop}{\ensuremath{\mathbf{Prop}}}
\newcommand*{\Set}{\ensuremath{\mathbf{Set}}}
\newcommand*{\sSet}{\ensuremath{\mathbf{sSet}}}

\newcommand*{\cgTop}{\ensuremath{\mathbf{Top}_{\mathrm{cg}}}}
\newcommand*{\blank}{\mathord{\hspace{1pt}\text{--}\hspace{1pt}}}

\mathchardef\texthyphen="2D
\newcommand*{\trType}[1]{\ensuremath{#1\mathbf{\texthyphen Type}}}
\newcommand*{\inftyoneCat}{\ensuremath{(\infty,1)\mathbf{\texthyphen Cat}}}
\newcommand*{\op}{^{\mathrm{op}}}
\newcommand*{\seq}{\overset{\mathrm{s}}{=}}
\newcommand*{\HoTTua}{HoTT$_{\mathrm{ua}}$}
\newcommand*{\HoTTel}{HoTT$_{\mathrm{el}}$}
\newcommand*{\HoTTelp}{HoTT$_{\mathrm{el}{+}{+}}$}
\newcommand*{\HoTTuf}{HoTT$_{\mathrm{UF}}$}
\DeclareMathOperator{\pt}{pt}

\DeclareMathOperator{\hasDimen}{hasDimension}

\DeclareMathOperator{\inl}{left}
\DeclareMathOperator{\inr}{right}
\DeclareMathOperator{\cglue}{glue}
\DeclarePairedDelimiter\trunc{\lVert}{\rVert} 
\DeclarePairedDelimiter\tr{\lvert}{\rvert} 
\DeclarePairedDelimiter\abs{\lvert}{\rvert} 
\DeclarePairedDelimiter\sem{\llbracket}{\rrbracket} 

\begin{document}


\title*{Higher Structures in Homotopy Type Theory}
\titlerunning{Higher Structures in HoTT}
\author{Ulrik Buchholtz}
\institute{Ulrik Buchholtz \at Fachbereich Mathematik,
  TU Darmstadt, 
  Schlossgartenstra\ss e 7, D-64289 Darmstadt,
  \email{buchholtz@mathematik.tu-darmstadt.de}}

\maketitle

\vspace*{-2.5cm}

\abstract{%
  The intended model of the homotopy type theories used in Univalent
  Foundations is the $\infty$-category of homotopy types, also known
  as $\infty$-groupoids. The problem of \emph{higher structures} is
  that of constructing the homotopy types needed for mathematics,
  especially those that aren't sets. The current repertoire of
  constructions, including the usual type formers and higher inductive
  types, suffice for many but not all of these. We discuss the
  problematic cases, typically those involving an infinite hierarchy of
  coherence data such as semi-simplicial types, as well as the problem
  of developing the meta-theory of homotopy type theories in Univalent
  Foundations. We also discuss some proposed solutions.}

\section{Introduction}
\label{sec:introduction}

Homotopy type theory is at the same time a foundational endeavor, in
which the aim is to provide a new foundation for mathematics, and an
area of mathematics and logic, in which the aim is to provide tools
for the mathematical analysis of homotopical (higher dimensional)
structures. Let us call the former Univalent Foundations (UF) and the
latter Homotopy Type Theory (HoTT), understanding that HoTT
encompasses many different particular type theories.

In the present chapter we use the issue of \emph{higher structures}
as a lens with which to study both of these aims and their relations
to other foundational approaches.

To motivate the problem of higher structures, we need to recall that
the intended universe of UF, and the principal model of HoTT, is the
realm of $\infty$-groupoids, a homotopical kind of algebraic
structures that have elements, identifications, identifications
between identifications, etc.\ ad infinitum, and these identifications
behave sensibly in that we can invert them, compose them, and whisker
by them, but the expected laws only hold up to higher identifications.
Grothendieck's \emph{homotopy hypothesis} tells us that
$\infty$-groupoids are the same as \emph{homotopy types}, so we shall
use these terms interchangeably, with a slight preference for the
latter, as then homotopy type theories are both theories of homotopy types as well as
homotopical type theories.

A common misconception is that higher homotopy types only occur in, or
only are relevant to, homotopy theory. That is very far from the case,
as even the type of sets, as used in most of mathematical practice, is
a $1$-type. And higher structures now feature prominently in many
areas ranging from geometry, algebra, and number theory, to the
mathematics of quantum field theories in physics and concurrency in
computer science. An introduction to homotopy types and the homotopy
hypothesis is given in Sect.~\ref{sec:inftygroupoids}.

It is a key point of difference between UF and earlier approaches to
foundations inspired by category theory that the former takes
$\infty$-groupoids rather than various notions of higher categories to
be the basic objects of mathematics, from which the rest are obtained
by adding further structure. This insight was due to
Voevodsky\footnote{In~\cite{Voevodsky2014} he wrote: ``The greatest
  roadblock for me was the idea that categories are `sets in the next
  dimension.' I clearly recall the feeling of a breakthrough that I
  experienced when I understood that this idea is wrong. Categories
  are not `sets in the next dimension.' They are `partially ordered
  sets in the next dimension' and `sets in the next dimension' are
  groupoids.''} who remarked that many natural constructions are not
functorial in the sense of category theory. (Think for example of the
center of a group.) However, every construction---if it is to have
mathematical meaning---has to preserve the relevant notion of
equivalence. (Isomorphic groups do have isomorphic centers, etc.)

Because UF aims to be a foundation for all of mathematics, it is
necessary that its language, in the shape of the HoTT, provide the means of
construction for all the homotopy types that are used in
mathematics. For the construction of \emph{sets}, this is not such a big
problem, as most of the sets that occur in mathematics can be
constructed from the type formers of Martin-L{\"o}f type theory. (But
even here there are subtleties if we wish to remain in the
constructive and predicative realm.)

The main problems appear when it comes to \emph{higher dimensional}
homotopy types. We discuss some positive results (structures that have
already been constructed) as well as some open problems (structures
that have not already been constructed) in Sect.~\ref{sec:problems}.

We remark that although we expect some actual \emph{negative} results
(i.e., impossibility proofs) for some of the open problems with
respect to some particular homotopy type theories, these have yet to appear. But anticipating that
further means of construction will be necessary, we discuss potential
solutions in Sect.~\ref{sec:solutions}.

For the remainder of this Introduction, we shall consider an analogy.
Martin-L{\"o}f type theory can be considered as a formal system for
making constructions. In fact, a variant with an impredicative
universe was called the \emph{Calculus of Constructions} (CoC), and a
further extension, the \emph{Calculus of Inductive Constructions}
(CIC) is the basis for the proof-assistant Coq.  And we shall be
concerned with the question of the limits of the methods of
construction available in constructive type theories.  An obvious
analogy presents itself, namely with euclidean geometry and the limits
of the methods of geometric constructions using ruler and
compasses. We shall (probably) find that, just as in the geometric
case, certain objects are not constructible from the most basic
constructions, and require further tools, such as the \emph{neusis},
for their construction. However, we shall follow Pappus'
prescription of parsimony and demand that everything that can be
constructed with lesser means, should be so constructed. As a
corollary, since a proof is a special case of a construction, we
demand that if something can be proved in a weaker system, then it
should be so proved.

Obviously we can include among the list of further means of
construction such well-known principles as the law of excluded middle
(LEM), Markov's principle (MP), the axiom of choice (AC), various kinds of transfinite
induction (TI), as well as principles of impredicativity. Some of
these, as well as weaker versions of these, are referred to as
\emph{constructive taboos} because admitting them is contrary to
certain philosophical outlooks inspired by constructivism or
intuitionism, and also because they cannot be mechanically executed at
all, or only with greatly increased computational complexity.

A further aspect of the constructive taboos is that they reduce the
number of \emph{models} in which we can interpret the
constructions. It is well-known that constructive systems admit many
useful models, indeed, this is one reason why classical mathematicians
may be interested in such systems. Non-homotopical constructive
systems can often be modeled in \emph{toposes}, more precisely,
$1$-toposes, which can be seen either as generalizations of Kripke
models, as generalized spaces, or indeed as generalized worlds of
sets.

It is suspected that HoTT can be modeled in higher toposes, more
precisely, \emph{\inftyone-toposes}. These dramatically extend the
usefulness of HoTT, for instance as explained in Schreiber's
Chapter. Earlier extensions of Martin-L{\"o}f type theory often
imposed axioms, such as the \emph{uniqueness of identity proofs}
(UIP), that rule out higher dimensional models. These contradict the
univalence axiom and may be called \emph{homotopical taboos}. More
refined axioms may hold in \inftyone-toposes corresponding to
$1$-toposes (the $1$-localic
\inftyone-toposes~\cite[Sect.~6.4]{LurieHTT}), but not in more general
\inftyone-toposes. These are called \emph{constructive-homotopical
  taboos}.

\section{Infinity groupoids and the homotopy hypothesis}
\label{sec:inftygroupoids}

The types in UF are supposed to be \emph{homotopy types}, so let us
dwell a bit on what they are, both from an intuitive point of view, and
from the perspective of mathematics developed in set-theoretic
foundations.

Intuitions are always hard to convey, and in the case of the notion of
homotopy type, even more so. Intuition is, after all, best developed
through practice and familiarity. One way to build an intuition for
homotopy types is through working in a homotopy type theory, either on
paper or with the help of a proof assistant. Many young workers in
HoTT/UF did this before learning about homotopy theory from a classical
point of view.

As a first approximation we can say that types $A$ are collections of
objects together with for each pair of objects $a,b:A$, a type of
\emph{identifications} $p : a=_A b$, together with meaningful
operations on these identifications, such as the ability to compose
and invert them. And there should also be higher order operations that
produce identifications between identifications, such as an
identification $\alpha(p) : (p^{-1})^{-1} = p$ for any $p : a =
b$. This description is meant to capture types in their incarnation as
\emph{$\infty$-groupoids}, and on this view, two types $A,B$ can be
identified if there is a (weak) functor $F:A \to B$ that is an
equivalence of $\infty$-groupoids.

Another intuition comes from describing types as (nice) topological
spaces up to homotopy equivalence. The objects are the points of the
space, and the identifications are the paths between points. 

The \emph{homotopy hypothesis} is the idea that these separate
intuitions capture the same underlying concept. It grew out of
Grothendieck's homotopy hypothesis concerning a particular
definition of $\infty$-groupoids~\cite{Grothendieck1983}. The modern
terminology is due to Baez~\cite{Baez2007}.

In order to explain the subtlety of the situation, let us turn to the
most common implementation of the idea of $\infty$-groupoids in the
context of set-theoretic mathematics. Here these are represented by
\emph{simplicial sets} satisfying a certain filling condition. These
simplicial sets are called \emph{Kan complexes} in honor
of~\cite{KanIII}. A simplicial set is a functor $X :
\Delta\op \to \Set$, where $\Delta$ is the category of
non-empty finite ordinals and order-preserving functions. This means
concretely that a simplicial set consists of a set of $n$-simplices $X_n$ for
each $n=0,1,\dots$ together with face and degeneracy maps satisfying
laws called the \emph{simplicial identities}. We think of the $0$
simplices as points, the $1$-simplices as lines between points,
$2$-simplices as triangles, etc.

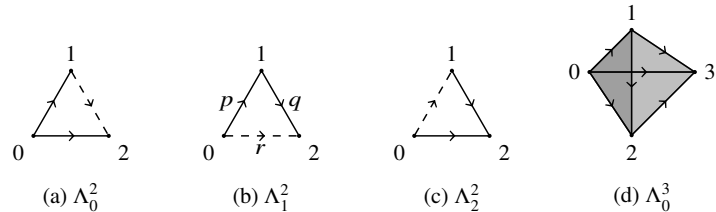
\begin{figure}
  \centering
  \begin{subfigure}[b]{0.21\linewidth}
    \centering
    \begin{tikzpicture}[dot/.style={draw,shape=circle,
        minimum size=1pt,inner sep=0pt,outer sep=0pt,fill=black}]
      \coordinate[dot,label=210:$0$] (p0) at (0,0);
      \coordinate[dot,label=$1$] (p1) at (60:1);
      \coordinate[dot,label=-30:$2$] (p2) at (1,0);
      \begin{scope}[decoration={markings,
          mark=at position 0.55 with {\arrow{>}}}]
        \draw[postaction={decorate}] (p0) -- (p1);
        \draw[postaction={decorate}] (p0) -- (p2);
        \draw[dashed,postaction={decorate}] (p1) -- (p2);
      \end{scope}
    \end{tikzpicture}
    \caption{$\Lambda^2_0$}    
  \end{subfigure}
  \begin{subfigure}[b]{0.21\linewidth}
    \centering
    \begin{tikzpicture}[dot/.style={draw,shape=circle,
        minimum size=1pt,inner sep=0pt,outer sep=0pt,fill=black}]
      \coordinate[dot,label=210:$0$] (p0) at (0,0);
      \coordinate[dot,label=$1$] (p1) at (60:1);
      \coordinate[dot,label=-30:$2$] (p2) at (1,0);
      \begin{scope}[decoration={markings,
          mark=at position 0.55 with {\arrow{>}}}]
        \draw[postaction={decorate}] (p0) -- (p1) node[midway,left] {$p$};
        \draw[dashed,postaction={decorate}] (p0) -- (p2)
          node[midway,below] {$r$};
        \draw[postaction={decorate}] (p1) -- (p2) node[midway,right] {$q$};
      \end{scope}
    \end{tikzpicture}
    \caption{$\Lambda^2_1$}\label{fig:KanComplex-b}
  \end{subfigure}
  \begin{subfigure}[b]{0.21\linewidth}
    \centering
    \begin{tikzpicture}[dot/.style={draw,shape=circle,
        minimum size=1pt,inner sep=0pt,outer sep=0pt,fill=black}]
      \coordinate[dot,label=210:$0$] (p0) at (0,0);
      \coordinate[dot,label=$1$] (p1) at (60:1);
      \coordinate[dot,label=-30:$2$] (p2) at (1,0);
      \begin{scope}[decoration={markings,
          mark=at position 0.55 with {\arrow{>}}}]
        \draw[dashed,postaction={decorate}] (p0) -- (p1);
        \draw[postaction={decorate}] (p0) -- (p2);
        \draw[postaction={decorate}] (p1) -- (p2);
      \end{scope}
    \end{tikzpicture}
    \caption{$\Lambda^2_2$}
  \end{subfigure}
  \begin{subfigure}[b]{0.21\linewidth}
    \centering
    \begin{tikzpicture}[dot/.style={draw,shape=circle,
        minimum size=1pt,inner sep=0pt,outer sep=0pt,fill=black},
      line join = round, line cap = round,
      x={(0:1cm)},y={(90:1cm)},z={(225:2mm)},scale=0.35]
      \pgfmathsetmacro{\factor}{1/sqrt(2)};
      \coordinate[dot,label=left:$0$] (p0) at (-2,0,-2*\factor);
      \coordinate[dot,label=above:$1$] (p1) at (0,2,2*\factor);
      \coordinate[dot,label=below:$2$] (p2) at (0,-2,2*\factor);
      \coordinate[dot,label=right:$3$] (p3) at (2,0,-2*\factor);
      
      \draw[fill=gray, opacity=.5] (p0.center)--(p1.center)--(p2.center)--cycle;
      \draw[fill=gray, opacity=.5] (p0.center)--(p1.center)--(p3.center)--cycle;
      \draw[fill=gray, opacity=.5] (p0.center)--(p2.center)--(p3.center)--cycle;
      \begin{scope}[decoration={markings,
          mark=at position 0.55 with {\arrow{>}}}]
        \draw[postaction={decorate}] (p0) -- (p1);
        \draw[postaction={decorate}] (p0) -- (p2);
        \draw[postaction={decorate}] (p0) -- (p3);
        \draw[postaction={decorate}] (p1) -- (p2);
        \draw[postaction={decorate}] (p1) -- (p3);
        \draw[postaction={decorate}] (p2) -- (p3);
      \end{scope}
    \end{tikzpicture}
    \caption{$\Lambda^3_0$}
  \end{subfigure}
  \caption{The Kan filling condition in dimensions 2 and 3.}
  \label{fig:KanComplex}
\end{figure}

The Kan filling condition says that if we are given $n$ compatible
$(n-1)$-simplices in $X$ in the sense that they could be $n$ of the
$n+1$ faces of an $n$-simplex, then there exists some such
$n$-simplex. This condition is illustrated in
Fig.~\ref{fig:KanComplex} in some low-dimensional cases. In each case,
we can think of the given data as a map from a \emph{horn}, a
sub-simplicial set $\Lambda^n_k\subseteq\Delta^n$ of the standard
$n$-simplex $\Delta^n$ consisting of the union of all the faces
opposite the $k$th vertex, into $X$. A lift is some extension of this
to a map from $\Delta^n$ to $X$, or equivalently, an $n$-simplex in
$X$ with the requisite faces.

For example, in Fig.~\ref{fig:KanComplex-b}, if we are given two
$1$-simplices $p$ and $q$ in $X$ with a common endpoint, then there
exists some $2$-simplex representing both a composite of $p$ and $q$
(the third face $r$)
together with the interior representing the fact that $r$ is the
composite of $p$ and $q$.

Note that Kan complexes give a \emph{non-algebraic} notion of
$\infty$-groupoid: there exists composites and higher simplicial
identifications, but there are no operations singling out a particular
composite.

Here we come to a potential pitfall: we cannot say that homotopy types
\emph{are} Kan complexes, for they have different criteria of
identity: In usual mathematical practice we identify two simplicial
sets if they are isomorphic (this is already a weaker notion of
identity than that provided by set theory!), whereas an identification
between Kan complexes $X$ and $Y$ \emph{considered as homotopy types}
should be a \emph{homotopy equivalence}. 

And this is perhaps an appropriate point at which to give a
type-theoretic take on Quine's~\cite{Quine1969} famous
slogans:\footnote{The second has also been discussed from a univalent
  perspective in \cite{Rodin2017,Tsementzis2017}.}
\begin{enumerate}
\item\label{it:Q1} \emph{To be is to be the value of a variable,} and
\item\label{it:Q2} \emph{No entity without identity.}
\end{enumerate}
In~\ref{it:Q1} we require moreover that all variables be typed, so we
say rather that to be an $A$ is to be the value of a variable of type
$A$ and more importantly, \emph{to be is to be an element of a type},
and in~\ref{it:Q2} we do not require any notion of identity between
entities of different types, but we do require as an essential part of
giving a type $A$ that the identity type $x=_Ay$, for $x,y:A$, is
meaningful and correctly expresses the means of identifying elements
of $A$: \emph{no type without an identity type.}

The discrepancy between the notion of identity between the model
objects (here Kan complexes) and the desired notion of identity (here
homotopy equivalence) is usually addressed using \emph{relative
  categories} as a tool. A relative category consists of a category
equipped with a wide subcategory of \emph{weak equivalences}. This is
often refined by adding more properties (e.g., the weak equivalences
satisfy the two-out-of-three or the two-out-of-six properties) or
structure, such as fibrations and/or cofibrations interacting nicely
with the weak equivalences. A particularly well-behaved notion is that
of a \emph{Quillen model category}, which does indeed contain both
fibrations and cofibrations in addition to weak equivalences, and is
assumed to be complete and cocomplete.

The category of simplicial sets $\sSet$ can be equipped with the
structure of a Quillen model category in which the fibrant objects are
the Kan complexes (these are also cofibrant as all objects are
cofibrant) and the weak equivalences between Kan complexes are the
homotopy equivalences.  The category of topological spaces, or more
precisely, for technical reasons, the category of compactly generated
topological spaces, $\cgTop$, can likewise be equipped with a Quillen
model structure in which the cofibrant objects are the nice spaces
(technically, cell complexes; all objects are fibrant) and the weak
equivalences between the nice spaces are the homotopy equivalences.

Quillen~\cite{Quillen1967} proved that these two model categories give
rise to equivalent \emph{homotopy categories}. For this purpose he
introduced the notion of (what is now called) a \emph{Quillen equivalence} between model
categories. Given a nice space $X$, the corresponding singular Kan
complex $\Pi_\infty(X)$ has as $n$-simplices the continuous maps from
the topological $n$-simplex into $X$, and given a Kan complex $A$, the
corresponding space is the \emph{geometric realization} $\abs A$ given
by gluing together topological simplices according to the face and
degeneracy maps in $A$.

That the homotopy categories are equivalent is a first step towards
getting what we actually want. We would actually like to show that the
Quillen model categories of simplicial sets and topological spaces
give rise to equivalent \emph{homotopy types} (in both cases
restricting to the objects that lie in a fixed Grothendieck
universe). And which notion of homotopy type should we use here? It
turns out not to matter, but it is easiest to make (large) Kan
complexes out of either one.

The way this is achieved is by enhancing both $\sSet$ and $\cgTop$ to
simplicially enriched categories (the latter via the singular Kan
complex construction on the level of mapping spaces) such that they
become \emph{simplicial model categories}, and then taking the
homotopy coherent nerves of the subcategories of homotopy equivalences
between bifibrant objects, i.e., between the model objects on both
sides.

Notice that to get a good theory of homotopy types in the classical
set-up we seem to need also a good theory of \inftyone-categories,
that is, categories (weakly) enriched in homotopy types, in order to
also get a good hold on the \emph{universe} of homotopy types, which
is another name of course anticipating the type-theoretic notion of a universe, and
which consists of the homotopy types that are small relative to some
Grothendieck universe.

There is another model structure on simplicial sets whose bifibrant
objects are the \emph{quasi-categories}, those that satisfy a weakening of
the Kan filling conditions that make them suitable as models of
\inftyone-categories. This notion was introduced by Boardman and
Vogt~\cite{BoardmanVogt1973} and the resulting theory of
\inftyone-categories has been studied extensively by
Joyal~\cite{Joyal2002} and Lurie~\cite{LurieHTT} (see also the
appendix of~\cite{LurieHTT} for details on simplicial model categories
as discussed above).

My point in bringing out these technicalities is not only to explain
how homotopy types are defined and handled in set-theoretic
mathematics, but also to give a sense of the subtleties involved. It
has taken many years to give a good account of how to treat higher
structure in set-theoretic mathematics (often by working in a
$1$-category-theoretic layer), and there are still many open questions
about which constructions and properties are invariant under weak
equivalences inside a model category and under Quillen equivalences
between model categories. For instance, it was just recently
established that a Quillen adjunction always induces an adjunction
between underlying quasi-categories, and hence an adjunction of the
presented \inftyone-categories~\cite{MazelGee2016}. Another line of
open questions concern the possibility of \emph{algebraic} models for
$\infty$-groupoids, where composition, inverses, etc., are given by
operations rather than merely assumed to exist. It is quite possible
that type theory will be influential in this area, see for instance
the suggestion of Brunerie~\cite[Appendix]{Brunerie2016}.

Thus it should come as no surprise that there are still open questions
about how to treat higher structures in HoTT/UF, which is a much
younger endeavor. These are the matters we shall now turn to.


\section{Higher Structures in HoTT/UF}
\label{sec:problems}

When Voevodsky proposed using type theory as a foundation for
mathematics, he based this on the insight that higher structures in
mathematics are \emph{not} always naturally objects of a higher
\emph{category}, but they \emph{are} always naturally objects of a
higher \emph{groupoid}.

Among the $\infty$-groupoids we find \emph{truncated} higher
groupoids, those whose structure is concentrated in a finite range of
dimensions. At the lowest level (truncation level $-2$) we find the
contractible types, those that only have one element up to
identifications. Secondly, we have propositions. These are types all
of whose identity types are contractible.

Moving up in the dimensions, we find next the sets, all of whose
identity types are propositions, and the $1$-groupoids, all of whose
identity types are sets, and so on. We recall from Altenkirch's
Chapter that these truncation levels have a natural formalization in
HoTT in terms of predicates
\[
  \hasDimen : \Ntwo \to \Type \to \Prop,
\]
and that we have corresponding types of $n$-truncated types, \trType n.

Not all types are truncated. The $2$-sphere, for example, has
structure in all dimensions, so it's not an $n$-type, for any $n$.

The $n$-types are related to the universe of all types, $\Type$, via
the truncation construction that maps a type $X$ to its closest
$n$-type $\trunc X_n$. There is a construction $\tr{\blank}_n : X \to
\trunc X_n$ giving rise to an equivalence
\[
  \blank \circ \tr{\blank}_n : (\trunc X_n \to Y) \to (X \to Y)
\]
for any $n$-type $Y$.

When we go to discuss higher structure, it is often the untruncated
types that are the hardest to construct. The principal reason
is that we can often construct truncated types in a top-to-bottom
fashion, dimensionwise. To construct a proposition, we can just
specify the type of evidence $P$ that the proposition is true and then if
necessary take the propositional truncation $\trunc{P}_{-1}$.

I want to emphasize at this point that the sets we discussed above
(and in the Chapters of Altenkirch and Ahrens-North) are \emph{not}
the sets of set theory! Following Quine's dictum, these are different
notions because they have different notions of identity. Let us
temporarily use subscripts to differentiate, and write set$_1$ for a
set theorist's set and set$_2$ for a structuralist/homotopy theorist's
set (this is also the model theorist's notion). There is even a third
notion of set, set$_0$, which is arguably more fundamental than either
set$_1$ or set$_2$, and which is the one taught in elementary
education.

From a type-theoretic point of view, a set$_0$ is simply a subset of a
fixed universal set$_2$, $X$. That is, we have the type
$\Set_0(X) :\equiv \mathcal P(X) :\equiv (X \to \Prop)$
representing the powerset of $X$. We have an elementary
membership relation, ${\in}_0: X \times \Set_0(X) \to \Prop$,
and two sets$_0$ are equal if they have the same elements in this
sense.

This is of course not the set-theorist's notion of set, according
which sets$_1$ are elements (rather than subsets) of a universe of
discourse $U$ (itself a set$_2$) that is equipped with a membership
relation ${\in}_1 : U \times U \to \Prop$ satisfying the axiom of
extensionality (and preferably many other set-theoretic axioms).

The naive set-theoretical hope would be to solve the equation $U =
\mathcal P(U)$ (as an identification of sets$_2$, i.e., an isomorphism).
This is impossible because of Cantor's diagonal argument, but it can
be approximated by the cumulative hierarchy $V$, a construction that can
be performed in HoTT via a higher inductive
type~\cite[Sect.~10.5]{HoTTBook}. Here $V$ is a large set$_2$ that is
the least solution of the equation
$V = \mathcal P_{\mathrm{small}}(V)$, where $\mathcal P_{\mathrm{small}}(V)
:\equiv \Sigma A : \Set. \Sigma f : A \to V. \mathrm{isInjective}(f)$
is the type of small subsets of $V$. Such sets$_1$ can be thought of
as certain well-founded trees, and their study has a quite
combinatorial flavor.

The default notion of set in HoTT/UF is set$_2$ given by the $1$-type
$\Set$, and this seems to be the one most often used in mathematical
practice outside of set theory. For instance, in almost all
mathematical contexts, each set can be replaced by an isomorphic copy
without changing the meaning of anything. Of course, sets$_0$ as
elements of powerset $0$-types/sets$_2$ also occur throughout
mathematical practice, but for these, the set theorist's and the
structuralist's notion coincide.

Likewise the notion of category splits into several distinct notions:
I will denote by \emph{precategory} the notion defined in Sect.~4.4 of
the Chapter by Ahrens-North, and leave the unadorned term
\emph{category} for a univalent precategory. Indeed, in most
category-theoretic contexts, each category can be replaced by an
equivalent while preserving the meaning. It is also useful to have
the term \emph{strict category}~\cite[Sect.~9.6]{HoTTBook} for a
precategory whose type of objects is a set. From the perspective of
set-theoretic mathematics, the $2$-type of categories arise from a
Quillen model category structure on the $1$-category of strict
categories.

Most of category theory can be formalized in HoTT/UF using the
univalent definition of category. A precategory can be thought of as a
category with extra structure, namely equipped with a functor from an
$\infty$-groupoid. For a strict category, this functor has as domain a
$0$-dimensional homotopy type. In set-theoretic foundations, it will
automatically be the case that every category can be equipped with
such a strict structure, but in UF this is an extra assumption, indeed
a constructive-homotopical taboo.

\subsection{Analytic and synthetic aspects of HoTT/UF}

A foundational theory must be synthetic, in that it describes how to
construct and reasons with its fundamental objects in terms of
postulated rules. It couldn't be otherwise, for if it described the
``fundamental'' objects in terms of other, more fundamental, objects, and
derived its rules from the properties of those, it would hardly be
foundational.

Homotopy type theories are synthetic theories of $\infty$-groupoids. The approach is
deeply \emph{logical}, where we think of \emph{logic as invariant
  theory} as pioneered by Mautner~\cite{Mautner1946} and later
developed by Tarski in a 1966 lecture~\cite{Tarski1986}. Both Mautner
and Tarski were inspired by the approach to geometry given in Klein's
\emph{Erlangen Program}~\cite{Klein1872}. The idea is that the logical
notions are those that are invariant under that maximal notion of
symmetries of the universes of discourse. If the universe of discourse
is a set, then the corresponding symmetry group is the \emph{symmetric
  group} consisting of bijections of the set with itself, but if it is
a higher homotopy type, then it is the (higher) \emph{automorphism}
group consisting of all self-homotopy equivalences.

In analogy with the synthetic theories of various notions of geometry
(euclidean, affine, projective, etc.), homotopy type theories are synthetic theories of
homotopy types (and set theories are synthetic theories of sets$_1$), cf.~also~\cite{Awodey2018}.

The analytic aspect is that all the rest of mathematics, all
mathematical objects, their types, and their structure, needs to be
developed in terms of homotopy types. And a key criterion for success
of a formalized notion is that it satisfies what Ahrens-North call the
\emph{principle of equivalence}, and which I linked to Quine's dictum
above: that the identity type captures the intended notion of
identifications between the mathematical objects that we are modeling.

One novel aspect of doing this analytical work in HoTT is when
defining a structured object, it can be a challenge already to get the
correct \emph{carrier type}. In set-theoretic foundations, any carrier
set of the correct cardinality will do, but in HoTT we are more
discerning.

We do reap some benefits of this extra care. For instance, any
construction (which, remember, could be proving a proposition,
inhabiting a set, etc.) we perform on a generic category is guaranteed
to be invariant under equivalence of categories, and we can use the
rules of identity types to transport the construction along any
equivalence.

Compare this to the situation in set-theoretic foundations: there we
have to prove invariance under equivalence for any construction on
types of dimension greater than $0$. For sets, this is not
necessary, because if we are given a set for which the notion of
identification between the elements is given by an equivalence
relation, we can take the quotient. In this way, the situation in set
theory is marginally better than that of type theory pre-HoTT, where
the set-quotient construction was not generally available, leading to
what some practitioners have called ``setoid hell''. But in set
theory, the same problem arises for any mathematical type of dimension
greater than zero, so we may surmise that formalizations based on set
theory will run into ``higher groupoid hell''.

\subsection{Some constructions that are possible}

Let us finally take a look at some constructions that are possible in
HoTT. Many of these are already discussed in~\cite{HoTTBook};
references are provided in other cases. I will structure this
discussion according to the means of construction used. Firstly, there
are those that only use the basic constructions in Martin-L{\"o}f type
theory, namely $\Sigma$- and $\Pi$-types, identity types, universes,
as well as (finitary) inductive types such as the natural numbers,
disjoint unions, the empty type, and the unit type.

Next, there are those that use in addition the univalence
axiom. Following that, there are those that can be reduced to one
particular higher inductive type, the (homotopy) pushout.

Finally, we find those constructions that seem to require more
advanced higher inductive types, and in the next subsection I shall
discuss those for which there is no known construction at the time of
writing.

In \emph{Basic Martin-L{\"o}f type theory} (MLTT) we can already
define many important notions such as homotopy fibers and other
pullbacks, the predicate $\hasDimen : \Ntwo \to \Type \to \Type$ and
the types \trType n. We have the types of categories and
$\dagger$-categories (cf.~\cite[Sect.~9.7]{HoTTBook}), as well as many
other types of mathematical objects occurring outside of
homotopy theory. But we are severely limited in our ability to
construct inhabitants of these types, or prove properties about
them. For instance, we cannot prove that $\hasDimen$ is valued in
propositions (this requires function extensionality), nor can we prove
that anything is not a set (such as the type $\Set$ itself), since
there a models of MLTT in which every type is a set. We cannot
construct set quotients, and, perhaps most embarrassingly, we cannot
even define the logical operations of disjunction and existential
quantification, as these require propositional truncation!

With \emph{univalence} we get function extensionality (as shown by
Voevodsky), and we can now prove many structural properties. Besides
$\hasDimen$ landing in $\Prop$, we can prove that \trType n{} is an
$(n+1)$-type, and we get the equivalence principle for the types of
algebraic structures and for categories as mentioned in the Chapter by
Ahrens-North. (See also~\cite{Awodey2014}.) We can also prove that the
$n$th universe is not an $n$-type for any external natural number
$n$~\cite{KrausSattler2015}.

At this point we can explore an intermediate route: instead of adding
higher inductive types, we can assume the \emph{propositional resizing
  axiom}~\cite[Axiom~3.5.5]{HoTTBook}, stating that the inclusion map
$\Prop_i\to\Prop_{i+1}$ of propositions in the $i$th universe into the
propositions in the $(i+1)$st universe is an equivalence. This makes
the theory impredicative, and it allows us to mimic many
impredicative tricks known from (constructive) set theory. For
instance, the propositional truncation can be defined as $\trunc A
:\equiv \Pi P:\Prop. (X\to P)\to P$. (If we have the law of excluded
middle, then we can just define $\trunc A :\equiv \lnot\lnot A$.)

The (homotopy) \emph{pushout} type is a simple, but versatile example
of a higher inductive type. It generalizes the disjoint union. Its
inputs are three types $A$, $B$, and $C$, together with functions
$f : C \to A$ and $g : C\to B$. (Such a configuration is called a
\emph{span}.) The pushout is a new type $D :\equiv A \sqcup^C_{f,g} B$ (often
written $A \sqcup^C B$ if $f$ and $g$ can be deduced from the context)
together with \emph{injections} $\inl{} : A \to D$ and $\inr{} : B \to
D$ fitting together in a square
\[
  \begin{tikzcd}
    C\ar[r,"g"]\ar[d,"f"'] & B\ar[d,"\inr{}"] \\
    A\ar[r,"\inl{}"'] & D
  \end{tikzcd}
\]
whose commutativity is given by a constructor $\cglue : \Pi
x:C. \inl(f\,x) = \inr(g\,x)$. (See~\cite[Sect.~6.8]{HoTTBook} for the
elimination and computation rules.)

Now a quite remarkable phenomenon appears. Most of the higher
inductive types that are commonly used can be constructed just from 
pushouts and the other constructions in MLTT with univalent
universes. These include joins and suspensions (and therefore, spheres),
cofibers (and thus smash products), sequential colimits, the propositional
truncation~\cite{vanDoorn2016,Kraus2016} and all the higher
truncations~\cite{Rijke2017}, set quotients, and in fact, all
non-recursive HITs specified using point-, $1$-, and $2$-constructors
by a construction due to van~Doorn~\cite{DRB2017}.
We also get cell complexes~\cite{BuchholtzFavonia2018},
Eilenberg-MacLane spaces~\cite{LicataFinster2014}, and projective
spaces~\cite{BuchholtzRijke2017}, and so a lot of
algebraic topology can be developed on this basis, and even a theory
of $\infty$-groups~\cite{BDR2018} and spectra (and thus homology and
cohomology theory), culminating in a proof that
$\pi_4(S^3)=\Z/2\Z$~\cite{Brunerie2016}, and a
formalized proof of the Serre spectral sequence for
cohomology~\cite{Doorn2018}.

Another important construction enabled by pushouts is the Rezk
completion, which turns a precategory into the category it
represents. This can in fact be done using univalence
alone~\cite{Ahrens-Kapulkin-Shulman2015}, at the cost of going to a
larger universe. This can be avoided either using pushouts or using
the propositional resizing axiom.

Because so many things can be developed on the basis of univalence,
pushouts, and propositional resizing, Shulman suggested that we define
an \emph{elementary \inftyone-topos} to be a finitely complete and
cocomplete, locally cartesian closed \inftyone-category with a
subobject classifier and object classifiers~\cite{Shulman2017}.  Let
me correspondingly introduce the term \emph{elementary HoTT} for MLTT
with univalence, pushouts, and propositional resizing.

Finally, let me mention some of the known constructions that seem to
require more than the above means, but that can nonetheless be
effected via more general higher inductive types. First, there is the
cumulative hierarchy as mentioned above~\cite[Sect.~10.5]{HoTTBook},
the Cauchy-complete real numbers~\cite[Sect.~11.3]{HoTTBook}, as well
as the partiality monad~\cite{AltDanKra2017}. I won't say more about
these, since this Chapter is supposed to be about \emph{higher}
structure. For these it is more relevant to mention
\emph{localizations} at a family of
maps~\cite{Rijke-Shulman-Spitters2017}. For example, if we localize at
a family of maps of the form $P(a) \to 1$, for $a:A$, where each type
$P(a)$ is a proposition, then we obtain a inner model of type theory
in itself, in this case a \emph{topological localization}.

Of course, the number of things that have been constructed and proved
in HoTT grows every day, so undoubtedly I'll have left some out.  Many
of these constructions have already been formalized in proof
assistants for HoTT. In the next section we move to constructions that
we don't yet know how to perform (and which may perhaps require new
means of construction).

\subsection{Some constructions that seem impossible}
\label{sec:negative-problems}

Because proving propositions is in HoTT/UF a special case of making
constructions, any currently open problems count as constructions we
don't yet know how to perform.\footnote{See
  \url{https://ncatlab.org/homotopytypetheory/show/open+problems} for
  an up-to-date list of open problems.} But for some of these, it is
expected that the difficulty is not just that the construction is
tricky to perform with the currently available means of
construction, but rather that we conjecture that entirely new means of
construction will be necessary.

The prime example that I will focus on is that of
\inftyone-categories, and the related notions of (semi-)simplicial
types. Intuitively, an \inftyone-category $\calC$ consists of a type
of objects $\calC_0$, for every pair of objects $a,b:\calC_0$ a type
of morphisms $\calC_1(a,b)$, operations for identities and
composition, operations that witness the unit- and associative laws,
operations that witness higher laws that these must satisfy, and so on
\emph{ad infinitum}. The problem is to come up with a way of specifying
all these higher coherence operations in a single type.

The basic example of an \inftyone-category from the point of view of
type theory is the category of types $\calS$ that has as type of
objects the universe $\Type$, and as morphisms from $A$ to $B$ the
type of functions from $A$ to $B$, $\calS_1(A,B):\equiv A \to
B$. Here, the evident identity and composition operations satisfy all
the laws and higher laws \emph{definitionally}, so it ought to be
particularly easy to show that $\calS$ is an \inftyone-category, that
is, if we knew how to define the type of \inftyone-categories,
\inftyoneCat.

It is crucial for the success of UF that we have a working definition
and theory of \inftyone-categories. Even if they are not
\emph{fundamental} in the sense that everything is built out of them,
they are still \emph{fundamental} in the sense that they are a key
tool in the development of modern higher algebra, geometry, and
topology.

The problem of defining \inftyone-categories is equivalent to the
problem of defining the type of simplicial types, \sType. A
\emph{simplicial type} is just a functor $X : \Delta\op\to\calS$, so
if we know how to define \inftyoneCat\ and the type of functors
between any two $\calC,\calD:\inftyoneCat$, then we can define
\sType. On the other hand, \inftyoneCat\ itself can be defined as the
subtype of \sType~consisting of \emph{complete Segal types} (also
called~\emph{Rezk types})~\cite{RiehlShulman2017}.

The problem of defining \inftyoneCat~can also be reduced to another,
apparently simpler problem, namely that of defining the type of
semi-simplicial types, \ssType. A \emph{semi-simplicial type} is a
functor $X:\Delta_+\op\to\calS$, where $\Delta_+$ is the subcategory
of $\Delta$ with the same objects but only injective
functions. However, $\Delta_+$ is a \emph{direct category},
viz.~a $0$-truncated category where the relation ``$x$ has a
nonidentity arrow to $y$'' is a well-founded relation on the
\emph{set} of objects, so we can give a more direct description as
follows: A semi-simplicial type $X$ consists of:
\begin{itemize}
\item a type of $0$-simplices $X_0$, and
\item for every pair of $0$-simplices $a_0,a_1:X_0$, a
  type of $1$-simplices from $a_0$ to $a_1$, $X_1(a_0,a_1)$, and
\item for every triple of $0$-simplices $a_0,a_1,a_2:X_0$ and
  $1$-simplices $a_{01}:X_1(a_0,a_1)$, $a_{02}:X_1(a_0,a_2)$, and
  $a_{12}:X_1(a_1,a_2)$, a type of $2$-simplices with boundary
  $a_{01},a_{02},a_{12}$,
\item \emph{and so on \ldots}
\end{itemize}
Again we have the problem that it seems that infinitely much data is
needed, but here it seems more plausible that an inductive (or
coinductive) approach could work. Only, no-one has figured out how to
do it, and at an informal poll of HoTT-researchers in Warsaw in 2015 a
majority believed that it is impossible.

We can define \inftyone-categories in terms of semi-simplicial types,
as the complete semi-Segal types~\cite{CapriottiKraus2017}. This work
was inspired by analogous work in the classical
setting~\cite{Harpaz2015}. Thus, the problems of defining
\inftyone-categories, simplicial types, and semi-simplicial types are
equivalent, but currently just beyond reach.

It is interesting to contrast the case of \inftyone-categories with
that of $\infty$-groups: An \inftyone-category structure on a pointed,
connected type of objects is the same as an $\infty$-monoid (the type
of which we also don't know how to construct). But the type of
$\infty$-groups is simply that of pointed, connected types, with the
type of group elements being the identity type
$\Omega A:\equiv(\pt=_A\pt)$, with $\pt:A$ the designated point.

While the above problems concern ``large'' types, there are also
problems of higher structure concerning ``small'' types. For example,
the three-sphere as a type $S^3$ should carry the structure of an
$\infty$-group, because it is the homotopy type of the Lie group
$SU(2)$. Thus, we should be able to construct the homotopy type of the
classifying space $BSU(2)$ with $\Omega BSU(2)=S^3$, but so far we've
not been able to do so. (We have the H-space structure, which is a
first approximation~\cite{BuchholtzRijke2018}.) In this case, however,
we expect that no new means of construction are needed.

An obvious approach would be to construct in the usual way a
simplicial set whose homotopy type is $BSU(2)$. If we then had the
realization operation $\abs\cdot : \sSet \to \Type$ that turns a
simplicial set into the homotopy type it represents, then we'd be
done. But such a realization operation itself seems impossible to
construct in elementary HoTT!

As a final, important, but more open-ended problem, let me mention the
problem of developing the meta-theory of HoTT inside HoTT/UF. This has
two sides, one relatively easy, and one quite hard. The relatively easy
side is the syntactic one, but even here there are difficulties. We
can represent \emph{extrinsic} untyped syntax and corresponding
transformations familiar from compiler theory: picking a surface
syntax, lexing and parsing this syntax, and then type-checking it. The
result should be \emph{intrinsic} syntax containing only well-typed
elements. The intrinsic syntax can be modeled by \emph{quotient
  inductive types} (QIT)~\cite{AltenkirchKaposi2016}, already
mentioned in Altenkirch's Chapter. The main difficulty here is one of
software-engineering: how do we structure both the intrinsic and
extrinsic syntax, and the transformations between them, in sufficient
generality to cover all the kinds of type theory we are interested
in.

The more difficult side is the semantic one. We want to define
interpretations of the intrinsic syntax in inner models, first of all
the canonical model where syntactic types $\vdash A$ are mapped to
types $\sem A$, syntactic terms $\vdash a:A$ are mapped to terms
$\sem a : \sem A$, and so on. (For proof-theoretic reasons, we expect
to only be able to represent the interpretation locally, for instance
type theory with $n$ universes inside the $(n+1)$st universe, or using
stronger principles in the target type theory than are in the source
type theory.) Shulman has called this problem ``making HoTT eat
itself''~\cite{Shulman2014}, for which the term \emph{autophagy}
suggests itself.

The problem is that if we use the QIT intrinsic syntax, then
everything syntactic is a (homotopy) set and the elimination rule will
only allow us to eliminate into sets, whereas for the canonical model
we're eliminating into \Type. And if we try to formalize the intrinsic
syntax using a non-truncated HIT, then it seems we need infinitely
many layers of coherence (reminding us of our problems above with
semi-simplicial types and the homotopical realization of simplicial
sets).

\section{Possible further means of construction}
\label{sec:solutions}

Now that we have seen concretely both the range of constructions that
are currently possible in (elementary) HoTT, and some prominent
problems that seem out of reach, let us take stock.

The first conclusion is that we'd very much like to \emph{prove} that
the problem of semi-simplicial types cannot be solved in elementary
HoTT. But assuming that, the next conclusion is that elementary HoTT
is by itself too incomplete to serve
its foundational role as the basis for UF. Further means of
constructions need to be added. But which ones, and how do we decide
which to add?

In some sense the situation is analogous to the question of new axioms
for set theory, but there are two main differences: First, we want to
use type theory as a \emph{programming language} and that means that
for any proposed extension, we should say how the new constructions
\emph{compute} when combined with the other constructions of type
theory. The univalence axiom is a sore point in this regard, as it has
been a long-standing open problem to give it a computational
meaning. This is now close to being solved via various \emph{cubical
  type
  theories}~\cite{Angiuli-Harper-Wilson2017,Bezem-Coquand-Huber2014,CCHM2018},
but there remains a question of whether the corresponding model
structures on various categories of cubical sets model
$\infty$-groupoids (we know that the \emph{test model structures}
do~\cite{BuchholtzMorehouse2017}). It is still completely open whether
the propositional resizing axiom can be given computational meaning.

Secondly, we want to use HoTT also in other models than
$\infty$-groupoids. It is conjectured that elementary HoTT can be
interpreted in any \inftyone-topos: a left-exact localization of the
functor category $\calC\op\to\calS$ for a small \inftyone-category
$\calC$. (It would take us too far afield to give the exact
formulation and up-to-date status of this conjecture;
see~\cite{HoTTWikiTopos}.)

Some of the most interesting targets are given by \emph{cohesive
  \inftyone-toposes}, whose objects can be thought of as geometrically
structured $\infty$-groupoids. (See also Schreiber's Chapter.) For
example, in the cohesive \inftyone-topos of smooth $\infty$-groupoids,
\SmoothInftyGpd, we find all smooth manifolds among the $0$-truncated
objects.  And we certainly want to be able to be able to reason about
smooth \inftyone-categories using a HoTT interpreted in
\SmoothInftyGpd.  So it will not do to propose a construction of the
type \inftyoneCat~that can't be performed meaningfully in any
\inftyone-topos.

In contrast, for the problem of interpreting type theory in internal
models, including these cohesive \inftyone-toposes, we need not require
that the means of doing so themselves are available in arbitrary
models. It seems sufficient to be able to do this ``at the top
level''. But however we solve this problem, it should probably be with
computationally meaningful (constructive) means, so that we're able to
do proofs by reflection inside these models. (Of course not all models
of interest will be definable constructively.)

Summing up, we expect there to be a stratification of homotopy type theories,

\smallskip

\begin{tabular}{@{\textbullet\ }l@{\ :\ }l}
  \HoTTua  & MLTT plus the univalence axiom and propositional resizing, \\
  \HoTTel  & \HoTTua\ plus pushouts, \\
  \HoTTelp & \HoTTel\ plus constructions needed for \inftyoneCat, \\
  \HoTTuf  & \HoTTelp\ plus reflective constructions.
\end{tabular}

\smallskip

\noindent Here, \HoTTua~is the basis for Voevodsky's UniMath
formalization effort~\cite{Voevodsky2015,UniMath}. We know that
\HoTTel~is strictly stronger (in the sense of having fewer models; not
in the sense of proof-theoretic strength), because it is consistent with
\HoTTua~that the $n$th universe is an $n$-type, while if the universes
are closed under pushouts, then they are not truncated. Since we don't
have impossibility proofs regarding the constructions of
\inftyoneCat~nor for autophagy, it is still conceivable that we can
take \HoTTelp~and \HoTTuf~to be \HoTTel.

For each of these we can consider adding classical axioms such as the
law of the excluded middle (LEM) or the axiom of choice (AC). These
can be seen as constructions that \emph{we} don't know how to perform
in general, but that an omniscient being would be able to perform. I
will not here take any sides in the debate about whether mathematics is
better done with or without these principles---I'll only say that it
seems to me that a foundational theory should give its users the
choice. (And in any case, we want theories without these constructive
taboos for reasoning about other toposes of interest.)

We can also remove the resizing axiom to get (generalized) predicative
systems, and for both the predicative and impredicative systems we may
add various generalized inductive types to increase the
proof-theoretic strength if needed, while keeping the systems
constructive and without changing their class of \inftyone-topos
models. Thus, the stratification above is meant primarily to
distinguish the (intended) models of the theories, and not their
proof-theoretic strength, which is an orthogonal concern.

With all that in mind, let us discuss some possible further means of
construction that we might add either for \HoTTel{} or \HoTTuf.

\subsection{Simplicial type theory}

For any \inftyone-topos $\calC$, we can consider the simplicial
objects in $\calC$, i.e., functors $\Delta\op \to \calC$, and this is
again an \inftyone-topos. As mentioned above, we find therein a full
subcategory of \inftyone-categories relative to $\calC$. This is the
basis for the suggestion by Riehl and Shulman~\cite{RiehlShulman2017}
for a synthetic type theory for \inftyone-categories. In their type
theory, let's call it sHoTT, types are interpreted as \emph{simplicial
  types}, and they give definitions for \emph{Segal} and \emph{Rezk
  types} with the latter representing \inftyone-categories. They can
also define a type of \emph{discrete} simplicial types, representing
ordinary types/$\infty$-groupoids, but this type is not a Rezk type
and so not the representation of the \inftyone-category of ordinary
types, $\calS$. Indeed, it is not clear at this moment, whether this
is at all representable in their system.

Much work remains before we can judge how useful this type theory is
for reasoning about \inftyone-categories. But from a philosophical
point of view it cannot be satisfactory to view sHoTT as a
foundational theory, for instance playing the role of \HoTTuf. (And it
is of course not intended as such!) Because simplicial types are
begging to be analyzed as such: simplicial objects in a category of
types, and not to be taken as unanalyzed in themselves. We really want
to be able to \emph{define} simplicial types inside a theory where
types are the fundamental objects.

\subsection{Two-level type theories}

Another approach to solving the problem of defining simplicial types
is to have another layer above the univalent type theory in which to
reason about infinitary strict constructions, including
(semi-)simplicial types. One proposal is Voevodsky's Homotopy Type
System (HTS)~\cite{Voevodsky2013}. This system has a distinction
between fibrant and non-fibrant types, and it has two identity types:
the usual homotopical identity type that only eliminates into fibrant
types as well as a new non-fibrant strict equality type that satisfies
the reflection rule: if $e:a\seq_A a'$ is an inhabitant of the strict
equality type, then $a$ and $a'$ are definitionally equal. This rule
of course makes type-checking undecidable, necessitating a further
language of evidence for typing derivations.

Another proposal is the two-level type theory
(2LTT)~\cite{AltenkirchCapriottiKraus2016,AnnenkovCapriottiKraus2017}. This
is similar to HTS in that it distinguishes between fibrant and
non-fibrant types (the latter are called \emph{pretypes}, but instead
of the reflection rule for the strict equality type it adds the rule
for uniqueness of equality proofs (Streicher's K) and function
extensionality as an axiom. Thus, type checking is decidable, but the
function extensionality axiom breaks computation.

In order to define simplicial types in 2LTT an extra principle beyond
the basic set-up is needed. This can be the assumption that the
fibrant and non-fibrant natural numbers coincide, or the more technical
assumption that Reedy fibrant diagrams of fibrant types indexed by a
strict Reedy category have fibrant limits. The former limits the class
of available models severely, while the latter does not.

However useful the two-level type theories may turn out to be, they
also seem unsatisfactory from a foundational perspective. Because what
is a pretype? Pretypes can only be motivated via the models of HoTT as
described in set-theoretic mathematics, where they arise as the
objects of a model category presenting an \inftyone-topos. But they
are not preserved by an equivalence of \inftyone-toposes, so they
don't have a presentation-independent meaning. They seem to be
merely a tool of convenience.

\subsection{Computational type theories}

If we limit ourselves to one constructive model, then there is a
principled way of making sense of new constructions. This is via the
paradigm of computational type theories in the Nuprl
tradition~\cite{Nuprl-book}. Here we consider a particular model to
give \emph{meaning explanations} for the judgments of type theory. For
a certain notion of cubical sets, this has been done by Harper's
group~\cite{Angiuli-Harper-Wilson2017,AngiuliHarper2018}. The benefit
of the approach is that it guarantees that all constructions are
computationally meaningful and make sense in the model. The downside
is that it is tied to a particular model (though by being judicious
with which primitives are added to the computational language this
downside can be minimized) and that it also leads to a type theory
without decidable type checking, so a separate proof theory or
language of evidence is needed.

\subsection{Presentation axioms}

It may have perhaps occurred to some readers that the problems
discussed in Sect.~\ref{sec:problems} should be solved in the same way
that they are solved in homotopical mathematics based in set-theory,
namely by working with set-based presentations.

We already mentioned geometric realization, an operation that produces
the underlying homotopy type of a given simplicial set or topological
space. We could consider adding $\abs\cdot : \sSet \to \Type$ as a
basic construction and the axiom stating that $\abs\cdot$ is
surjective, meaning that every type is merely equivalent to the
geometric realization of some simplicial set. And perhaps we should
add a further axiom stating that every function $A \to B$ between
types arises (merely) as the geometric realization of a function
between representing simplicial sets.

Something like this may indeed be appropriate at the level of
\HoTTuf~if it could be given a computational meaning. But certainly
not at the level of \HoTTelp, because the axioms would severely
restrict the range of models (they are constructive-homotopical
taboos).

Even the much weaker axiom \emph{sets cover} (SC), stating that every
type admits a surjection from a set admits a simple counter-model
\inftyone-topos~\cite{nLabTypesCover}.

On the other hand, SC (or something like it) is necessary in order to
describe the semantics of HoTT (with universes) in presheaf
toposes. Indeed, the universe in a presheaf topos is built from
certain \emph{sets} covering the $1$-types of presheaves of small
sets.

\section{Conclusion}
\label{sec:conclusion}

Higher structures are at once the raison d'{\^e}tre and so far, the
Achilles' heel, of HoTT/UF from a foundational perspective. HoTT can
handle with ease many important higher structures, such as the
$1$-type of sets and the $2$-type of categories, that can only
imperfectly be represented in other foundational systems. But so far
it cannot define the (untruncated) type of \inftyone-categories, and
this is a major impediment to the foundational aspirations of
HoTT/UF. Of course, HoTT can be (and has been) used successfully to
reason about (structured) homotopy types. In this way, the various
homotopy type theories function as \emph{domain specific languages} (DSLs).

To be foundational, however, we need to find a compelling construction
of, and theory of, \inftyone-categories, and of the semantics of
HoTT-DSLs, inside homotopy type theory itself. It appears that new
methods of construction are needed, but it is at this time not clear
what they should be.

A dramatic possibility, not mentioned in Sect.~\ref{sec:solutions}, is
that we should take \inftyone-categories to be fundamental after all,
and build a synthetic type theory where the types are
\inftyone-categories rather than $\infty$-groupoids. This would be a
\emph{directed} type theory. Such a thing would undoubtedly be quite
complicated due to the need to keep track of variances
(see~\cite{LicataHarper2011,Nuyts2015} for some preliminary attempts),
and it would represent a return to the old ways of thinking about
categorical foundations, albeit updated to account for a homotopical
perspective. We would carve out the $\infty$-groupoids as those
\inftyone-categories all of whose morphisms are invertible rather than
trying to build \inftyone-categories out of $\infty$-groupoids. In any
case, directed type theories should be useful also as DSLs for
reasoning about $(\infty,2)$-toposes.

Personally, I think we'll find some solution that allows us to stay at
the level of $\infty$-groupoids for the foundational theory. Perhaps
there is a kind of two-level type theory that allows us to capture
the strict nature of the \inftyone-category of types without
postulating a bunch of meaningless pretypes.

An analogy can perhaps be made with the foundations of stable homotopy
theory. The \inftyone-category of spectra is a symmetric monoidal
stable \inftyone-category, and from a foundational point of view, this
is the correct viewpoint, since spectra should be identified when they
are weakly equivalent. However, it was discovered that this
\inftyone-category can be presented by symmetric monoidal Quillen
model categories (i.e., very strict structures), and this has been
very important in facilitating computations in stable homotopy
theory~\cite{EKMM1995,MMSS2001}. (I should mention that there is
work-in-progress by Finster-Licata-Morehouse-Riley on developing a
HoTT-DSL for stable homotopy theory targetting the cohesive
\inftyone-topos of parametrized spectra: this captures the strict
monoidal structure of spectra in a type theory.)

And so it may be, that in order to realize the foundational potential
of HoTT/UF, we shall need to capture the strict structure of type
theory itself, perhaps by reflecting more of judgmental structure at
the level of types.

I'm confident that a good solution will be eventually found. The field
is still young, and it will be exciting to see what the future brings.

\bibliographystyle{spmpsci}
\bibliography{fomus}
\end{document}